\UseRawInputEncoding
\documentclass[review,3p,10pt]{elsarticle}

\usepackage{natbib}

\usepackage{amssymb}
\usepackage{amsmath}
\usepackage{cases}
\usepackage{epstopdf}
\usepackage{fixltx2e}
\usepackage{mathtools}
\usepackage{bm}
\usepackage{multirow}
\usepackage{color}
\usepackage{algorithm}
\usepackage{algorithmic}
\usepackage{amsthm}
\usepackage{mathrsfs}
\usepackage{amssymb}

\newtheorem{remark}{Remark}
\usepackage{lineno}

\usepackage{url}

\newcommand{\tabincell}[2]{\renewcommand\arraystretch{0.8}\begin{tabular}{@{}#1@{}}#2\end{tabular}}

\journal{arXiv.org}

\begin{document}
\begin{frontmatter}

\title{Extended Load Flexibility of Utility-Scale P2H Plants: Optimal Production Scheduling Considering Dynamic Thermal and HTO Impurity Effects }

\author[label1]{Yiwei~Qiu}
\author[label1]{Buxiang~Zhou}
\author[label1]{Tianlei~Zang\corref{cor1}}
\ead{zangtianlei@scu.edu.cn}
\author[label1]{Yi~Zhou}
\author[label1]{Shi~Chen}
\author[label2]{Ruomei~Qi}
\author[label2]{Jiarong~Li}
\author[label2]{Jin~Lin}

\address[label1]{College of Electrical Engineering, Sichuan University, Chengdu, 610065, China}
\address[label2]{State Key Laboratory of Control and Simulation of Power Systems and Generation Equipment, Department of Electrical Engineering, Tsinghua University, Beijing 100087, China}
\cortext[cor1]{Corresponding author}

\begin{abstract}
  In the conversion toward a clear and sustainable energy system, the flexibility of power-to-hydrogen (P2H) production enables the admittance of volatile renewable energies on a utility scale and provides the connected electrical power system with ancillary services. To extend the load flexibility and thus improve the profitability of green hydrogen production, this paper presents an optimal production scheduling approach for utility-scale P2H plants composed of multiple alkaline electrolyzers. Unlike existing works, this work discards the conservative constant steady-state constraints and first leverages the dynamic thermal and hydrogen-to-oxygen (HTO) impurity crossover processes of electrolyzers. Doing this optimizes their effects on the loading range and energy conversion efficiency, therefore improving the load flexibility of P2H production. The proposed multiphysics-aware scheduling model is formulated as mixed-integer linear programming (MILP). It coordinates the electrolyzers' operation state transitions and load allocation subject to comprehensive thermodynamic and mass transfer constraints. A decomposition-based solution method, SDM-GS-ALM, is followingly adopted to address the scalability issue for scheduling large-scale P2H plants composed of tens of electrolyzers. With an experiment-verified dynamic electrolyzer model, case studies up to 22 electrolyzers show that the proposed method remarkably improves the hydrogen output and profit of P2H production powered by either solar or wind energy compared to the existing scheduling approach.
\end{abstract}

\begin{keyword}
    alkaline water electrolysis  \sep  demand side management  \sep  hydrogen production  \sep  augmented Lagrangian method  \sep  production scheduling  \sep  unit commitment
\end{keyword}

\end{frontmatter}





\section{Introduction}
\label{sec:intro}




\subsection{Motivation}
\label{sec:motivation}

Utility-scale hydrogen production via water electrolysis has been recognized as a promising path toward renewable energy admittance and the decarbonization of the chemical and transportation industries \cite{jovan2022utilization,
li2021co,klyapovskiy2021optimal}. Driven by the demand for \emph{green hydrogen substitution} in these now carbon emission-intensive sectors, the demand for green hydrogen is forecasted to exceed 40 Mt/yr in China alone in 2050 \cite{alliance2019white}. Demonstration projects of utility-scale renewable power-based hydrogen production and industrial utilization (such as ammonia and methanol synthesis) have been approved or under construction around the world \cite{Helios2022,Damao2022,Lanzhou2018will}.

As an electrical power load, power-to-hydrogen (P2H) production needs to accommodate the fluctuating wind or solar energy \cite{zhang2014wind,li2020capacity,serna2017predictive,fang2019control} or provide auxiliary services, e.g., peak shaving and frequency regulation, for the power system \cite{
kopp2017energiepark,el2019hydrogen,dozein2021fast}. For example, in Inner Mongolia, China, it is required that green hydrogen production consumes at least $80\%$ of local renewable energy \cite{inner2022}. The latest European Commission regulatory draft requires that hydrogen production fully uses renewable energy \cite{ec2022}. These policies underscore the importance of the load flexibility of P2H production.

The load flexibility, quantified by the loading range, ramping rate, and energy conversion efficiency in varying-load operation, reflects the P2H plant's ability to accommodate volatile renewable energy or provide ancillary regulatory services. As a result, better flexibility means higher hydrogen output and profitability \cite{serna2017predictive,fang2019control}.
Therefore, in production scheduling, the flexibility of the P2H plant needs to be fully exploited.

Water electrolysis is the most common method for renewable energy-based P2H production \cite{david2020dynamic}. Mainstream technical routes include alkaline water electrolysis (AEL), proton exchange membrane electrolysis (PEMEL), and solid oxide cell electrolysis (SOCEL) \cite{grigoriev2020current}. Due to its relatively high maturity, large capacity, and long lifespan, AEL is preferred by many for industry-scale green hydrogen production \cite{grigoriev2020current,varela2021modeling,straka2021comprehensive}. This work focuses on AEL-based P2H production.

\begin{figure}[t]
  \centering
  \includegraphics[width=6.2in]{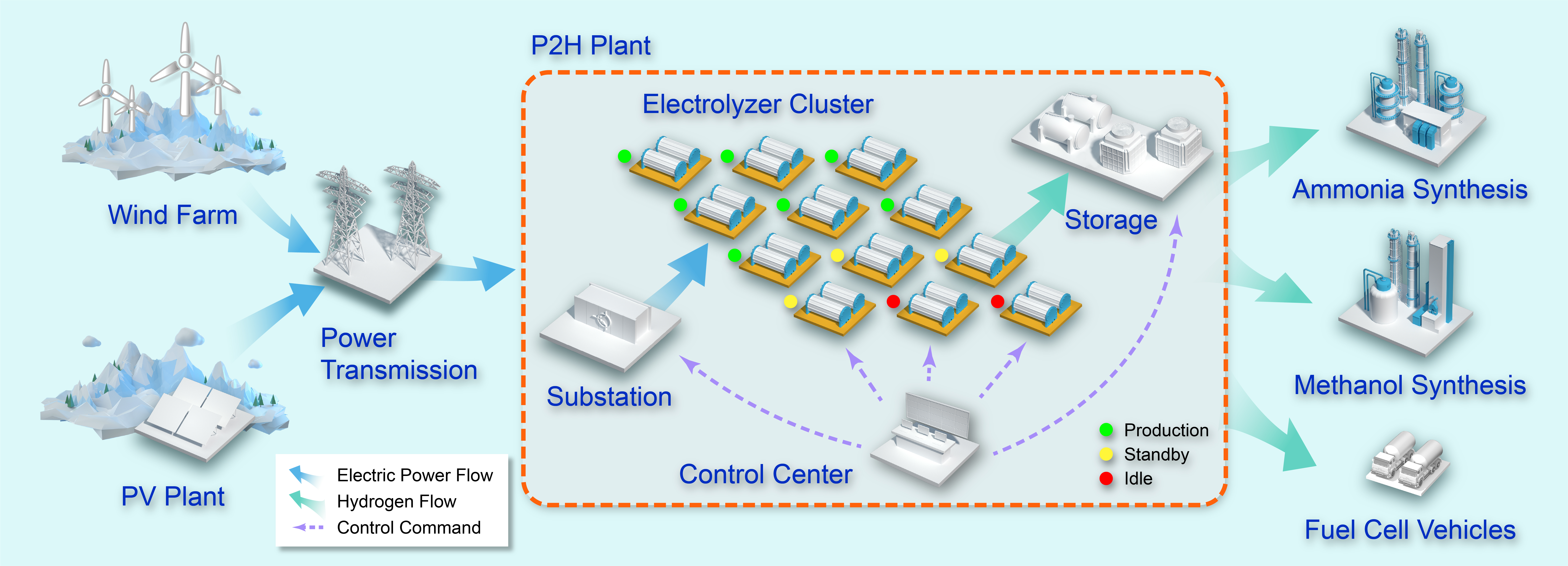}
  \caption{The conceptual schematic of utility-scale hydrogen production.}
  \label{fig:cluster}
\end{figure}

The conceptual schematic of an utility-scale renewable power-based hydrogen production system is shown in Fig. \ref{fig:cluster}. The load flexibility of P2H production is subject to two aspects. The first is the complicated dynamic thermal and mass transfer constraints \cite{david2020dynamic}. For example, to avoid excessive overvoltage, an alkaline electrolyzer cannot be fully loaded unless it has warmed up; it cannot operate at a low load level for a long duration due to the accumulation of hydrogen-to-oxygen (HTO) impurities that may cause a flammable mixture; and its energy efficiency is significantly affected by the temperature.
In addition, due to the limited capacity of a single electrolyzer, an utility-scale P2H plant is composed of many electrolyzers \cite{grigoriev2020current,varela2021modeling}. Hence, the following two aspects need to be determined in the production scheduling:
\begin{enumerate}
  \item{Planning the operational state transitions}: similar to the classical power system unit commitment (UC) problem, the states of the electrolyzers, i.e., Production, Idle, and Standby, need to be determined to fit the power supply profile.
  \item{Allocating the load}: the power supply needs to be allocated to each electrolyzer to maximize hydrogen production and ensure the multiphysics constraints are not violated.
\end{enumerate}

Although renewable power-based P2H production scheduling is a hot topic, most existing works use an overly simplified electrolyzer model, neglecting the dynamic multiphysics constraints. Unlike those works, this work first jointly considers the dynamic temperature and HTO impurity effects in scheduling. Compared to the traditional methods with a fixed load range and ramping limits, this work fully exploits the load flexibility of AEL-based P2H production. A brief literature review is presented in Section \ref{sec:review}, and the contribution of this work is given in Section \ref{sec:contribution}.

\begin{table*}[tb]\scriptsize
  \renewcommand{\arraystretch}{1.3}
  \caption{Summary of the Latest Literature of Power-to-Hydrogen Production Scheduling and Control}\vspace{6pt}
  \label{tab:literature}
  \centering
  \begin{tabular}{cccccccc}
  \hline \hline
  \multirow{2}{*}{Literature}                    & \multirow{2}{*}{\tabincell{c}{Number of \\ Electrolyzers}}
                                                 & \multirow{2}{*}{\tabincell{c}{Type}}    & \multicolumn{4}{c}{Electrolyzer Model}                                         & \multirow{2}{*}{\tabincell{c}{Solution\\Method}} \\
  \cline{4-7}                                    &                           &             & \tabincell{c}{State\\Transition}  & \tabincell{c}{Production\\Function} & \tabincell{c}{Temperature}       & \tabincell{c}{HTO\\  Impurity}              \\                        \hline
  Serna 2017 \cite{serna2017predictive}          & Multiple                  & AEL         & \checkmark        & Nonlinear           & $\times$          & $\times$         & MIQP   \\ 
  Fang 2019 \cite{fang2019control}               & Single                    & PEMEL       & \checkmark        & Nonlinear           & $\times$          & $\times$         & \tabincell{c}{Heuristic\\ (Rule-Based)}  \\
  Hong 2022 \cite{hong2022optimization}          & Multiple                  & AEL         & \checkmark        & Nonlinear           & $\times$          & $\times$         & Fuzzy Control \\
  Klyapovskiy 2021 \cite{klyapovskiy2021optimal} & Single                    & Unspecified & $\times$          & Linear              & $\times$          & $\times$         & MILP   \\
  Uchman 2021 \cite{uchman2021varying}           & Multiple                  & Unspecified & \checkmark        & Linear              & $\times$          & $\times$         & \tabincell{c}{Exhaustive\\Search} \\
  Varela 2021 \cite{varela2021modeling}          & Multiple                  & AEL         & \checkmark        & Linear              & $\times$          & $\times$         & MILP \\
  Li 2022  \cite{li2021co}                       & Single                    & AEL         & $\times$          & Linear              & $\times$          & $\times$         & MILP \\
  He 2021  \cite{he2021hydrogen}                 & Multiple                  & Unspecified & \checkmark        & Linear              & $\times$          & $\times$         & MILP \\
  Flamm 2021 \cite{flamm2021electrolyzer}        & Single                    & PEMEL       & $\times$          & Nonlinear           & \checkmark        & $\times$         & MILP    \\
  Zheng 2022 \cite{zheng2022optimal}             & Single                    & AEL         & \checkmark        & Linear              & \checkmark        & $\times$         & MILP  \\
  Shen 2021 \cite{shen2021coordination}          & Multiple                  & AEL         & \checkmark        & Linear              & \checkmark        & $\times$         & \tabincell{c}{Heuristic\\  (Rule-Based)}  \\
  Yang 2022 \cite{yang2022scheduling}            & Single                    &  AEL        & \checkmark        & Nonlinear           & $\times$          & $\times$         & \tabincell{c}{Mixed-Logic\\MILP}  \\  \hline
  \bf This work                          & \bf Multiple              & \bf AEL     & \bf \checkmark    & \bf Nonlinear       & \bf \checkmark    & \bf \checkmark   &  \bf \tabincell{c}{Decomposition-\\Based MILP} \\
  \hline \hline
  \end{tabular}
\end{table*}

\subsection{A Brief Literature Review}
\label{sec:review}

The community has long recognized the flexibility of P2H production. Researchers have carried out technical-economic analyses for renewable energy admittance \cite{shams2021machine,genovese2021parametric,bhandari2021hydrogen,fragiacomo2020technical} or providing ancillary services \cite{kopp2017energiepark}. However, the P2H plant is generally modeled as an energy node with a constant conversion factor. The detailed operational feasibility of the electrolyzers is omitted. Although such simplification is acceptable for technical-economic analysis, it does not provide information on how to operate a plant.

To investigate the operation of a P2H plant to accommodate volatile renewable energy, production scheduling approaches have been proposed. For example, Serna et al. \cite{serna2017predictive} proposed an energy management model for offshore hydrogen production powered by wind and wave energy. 
Fang et al. \cite{fang2019control} developed a rule-based strategy for coordinating the electrolyzer and energy storage. Hong et al. \cite{hong2022optimization} proposed a fuzzy controller for tracking maximal hydrogen production using wind power. However, its rule-based strategies cannot ensure an optimum. Klyapovskiy et al. \cite{klyapovskiy2021optimal} proposed an energy management framework for ammonia production using green hydrogen. It used an aggregated factor to model the relation between electricity and hydrogen.
Uchman et al. \cite{uchman2021varying} used an exhaustive search to determine the production plan of three electrolyzers. However, the solution method is difficult to scale up. Varela et al. \cite{varela2021modeling} used mixed-integer programming (MILP) to determine the on-off state transitions of multiple electrolyzers in a plant. 
Li et al. \cite{li2021co}, He et al. \cite{he2021hydrogen}, and Ahmadi et al., \cite{ahmadi2022design} further combined the scheduling of P2H with the storage, transportation, and utilization of hydrogen.

The classical literature \cite{ulleberg2003modeling} and engineering practices have pointed out that the flexibility and efficiency of an electrolyzer are affected by its thermodynamic and mass transfer constraints. However, most existing scheduling approaches omit these multiphysics dynamics. Instead, the electrolyzers are modeled with a constant loading range, ramping rate, and efficiency factor, leading to overly conservativeness of the scheduling result. For example, the electrolyzer can physically operate at a low load level for a short time without violating the impurity constraint, as explained in Section \ref{sec:hto}. Unfortunately, this feature cannot be utilized in the current scheduling models, as a constant lower loading limit is enforced.

Recently, to improve flexibility, varying-load control of alkaline electrolyzers considering multiphysics dynamics has drawn attention. Flamm et al. \cite{flamm2021electrolyzer} and Zheng et al. \cite{zheng2022optimal} developed model predictive controllers (MPCs) to address temperature dynamics and their effect on efficiency. Qi et al. \cite{qi2021pressure} proposed a pressure control strategy to alleviate HTO impurity accumulation to extend the lower loading limit. Yang et al. \cite{yang2022scheduling} used a detailed CFD-based lye-bubble two-phase flow model to obtain an accurate steady-state efficiency curve. However, these works focus on a single electrolyzer in small-scale applications. 

A review of the recent literature on P2H production scheduling is summarized in Table \ref{tab:literature}. As summarized, there is a lack of a scheduling approach for an industry-scale P2H plant to coordinate multiple electrolyzers in the plant while being able to leverage the electrolyzer thermal and mass transfer dynamics. Filling the gap to improve the flexibility and profitability of renewable energy-based hydrogen production is the target of this work.

\begin{figure}[t]
  \centering
  \includegraphics[scale=1.2]{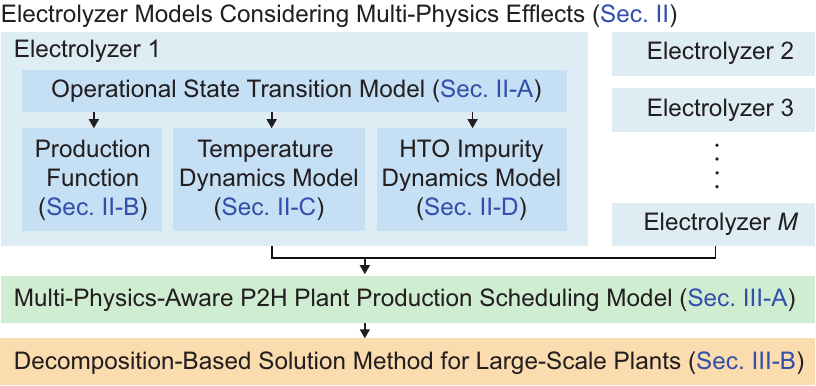}
  \caption{The conceptual structure of this paper.}
  \label{fig:schematic}
\end{figure}

\subsection{Contributions of This Work}
\label{sec:contribution}


This paper\footnote{This paper is a substantial extension of the 6-page conference paper \cite{qiu2022extended}.} aims to fully exploit the load flexibility of a renewable energy-powered P2H plant. According to the above discussions, the varying-load operation of an industry-scale P2H plant is subject to complicated multiphysics constraints, including temperature and HTO impurity effects. They have a nonneglectable impact on flexibility and efficiency. A comprehensive scheduling model is required to fit the varying load operations with a volatile power supply. In addition, an efficient solution method is needed to tackle the large-scale and complicated multiphysics-aware scheduling model. To address these requirements, an approach that considers the dynamic multiphysics feasibility requirement of electrolyzers is proposed. The main contributions of this work include the following:
\begin{enumerate}
    \item A comprehensive multiphysics-aware scheduling model is first presented for large-scale P2H plants. The dynamic thermal and mass transfer constraints are considered instead of the traditional fixed steady-state constraints. 

    \item 
        To address the scalability issue for scheduling a large P2H plant with multiple electrolyzers, a decomposition-based solution method, namely, SDM-GS-ALM \cite{boland2019parallelizable}, is adopted to leverage the natural separability of the scheduling problem and solve it effectively.

    \item Case studies show that for a plant directly coupled with wind or solar energy, the profit of hydrogen production is improved by $1.438\%$ and $0.982\%$ by considering the multiphysics dynamics.
\end{enumerate}


The structure of this paper is given in Fig. \ref{fig:schematic}. Section \ref{sec:constraint} presents the dynamic multiphysics constraints of an alkaline electrolyzer; Section \ref{sec:scheduling} presents the overall production scheduling model and the decomposition-based solution method; finally, Section \ref{sec:case} verifies the proposed scheduling approach by case studies.

\section{Modeling Dynamic Multiphysics Constraints of an Alkaline Electrolyzer for Scheduling}
\label{sec:constraint}

\begin{figure}[t]
  \centering
  \includegraphics[scale=1.15]{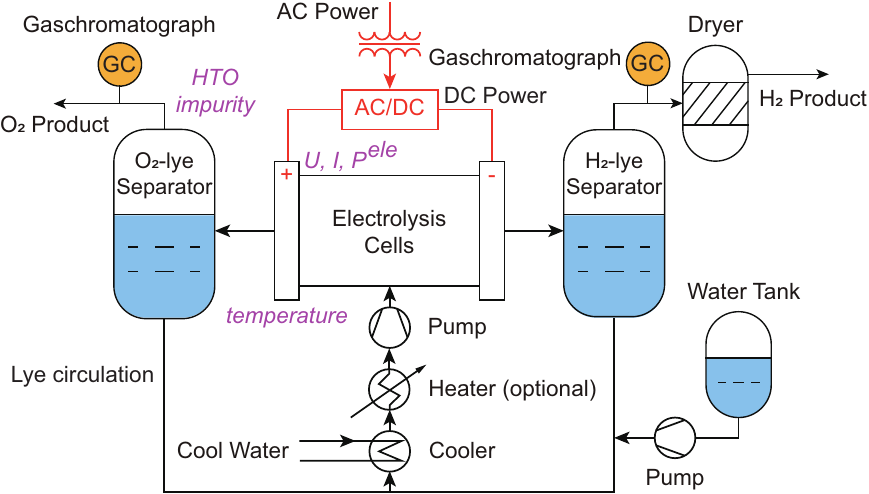}
  \caption{Illustrative schematic of an alkaline electrolyzer.}
  \label{fig:ale}
\end{figure}


As illustrated in Fig. \ref{fig:ale}, an alkaline electrolyzer comprises the electrolysis cells (stack) and auxiliary equipment, including lye-gas separators, heat exchangers, pumps, and power supply; an elaborated description can be seen in \cite{olivier2017low}. Moreover, an industry-scale P2H plant comprises multiple electrolyzers. For example, to meet the hydrogen demand of an ammonia plant rated $100$ kt/yr, the P2H plant needs $20$ electrolyzers, each with a rated hydrogen production of 1,000 $\text{Nm}^{2}/\text{h}$ (about $5$ MW rated power).

To fit the fluctuating renewable energy, the number of electrolyzers in production changes with time. The electrolyzers switch between three operational states, namely, \emph{Production (P)}, \emph{Standby (S)}, and \emph{Idle (I)} \cite{varela2021modeling,zheng2022optimal}. In \emph{Production}, the electrolyzer breaks up water molecules into hydrogen and oxygen, the pump keeps lye circulating, the cooler takes away excess heat, and the control system keeps the temperature and pressure at appropriate levels \cite{olivier2017low}. The total power consumption includes electrolytic power and auxiliary consumption. In \emph{Standby}, the electrolytic power is zero, but the auxiliary system keeps working so that it can quickly switch to \emph{Production}.
In \emph{Idle}, the system is switched off, and no power is consumed.

The state transition, thermodynamics, and mass transfer processes are coupled in scheduling. The temperature and HTO impurity remarkably impact the electrolyzer efficiency, ramping rate, and load range \cite{zheng2022optimal,flamm2021electrolyzer,qi2021pressure}. To construct a comprehensive plant-level scheduling model, Section \ref{sec:state} describes the state transition of the electrolyzers; Section \ref{sec:prod} establishes the relation between power consumption and hydrogen production; and the dynamic temperature and mass transfer models are given in Sections \ref{sec:temp} and \ref{sec:hto}, respectively.

We denote the scheduling horizon and step length as $N$ and $h$, respectively, and the number of electrolyzers as $M$. The subscripts $i,j$ indicate a quantity for the $i$th electrolyzer at time step $j$. The multiphysics feasibility constraints are given below.



\subsection{Production/Idle/Standby State Transition of Electrolyzers}
\label{sec:state}

Following Varela et al. \cite{varela2021modeling} and Zheng et al. \cite{zheng2022optimal}, the state of the $i$th electrolyzer at step $j$ is represented by three mutually exclusive binary variables $b_{i,j}^{\text{P}}$, $b_{i,j}^{\text{S}}$, and $b_{i,j}^{\text{I}}$, as
\begin{align}
  b_{i,j}^{\text{P}} + b_{i,j}^{\text{S}} + b_{i,j}^{\text{I}} = 1. \label{eq:states}
\end{align}


The \emph{Startup}, i.e., switching from \emph{Idle} to \emph{Production} or \emph{Standby}, is indicated by a binary variable $b_{i,j}^{\text{SU}}$, following
\begin{align}
   b_{i,j}^{\text{P}} +  b_{i,j}^{\text{S}} +  b_{i,j-1}^{\text{I}} - 1 &\le b_{i,j}^{\text{SU}}.  \label{eq:startup}
\end{align}

Meanwhile, a minimal gap of $N^{\text{gap}}$ steps between shutdown and startup is enforced by
\begin{align}
 b_{i,k-j}^{\text{I}} + b_{i,j}^{\text{I}} - \sum_{l=1}^{k-1}  b_{i,j-k+l}^{\text{I}} \le 0,\ \forall k = 2, \ldots,  N^{\text{gap}}. \label{eq:gap}
\end{align}

\begin{figure}[t]
  \centering
  \includegraphics[scale=1.15]{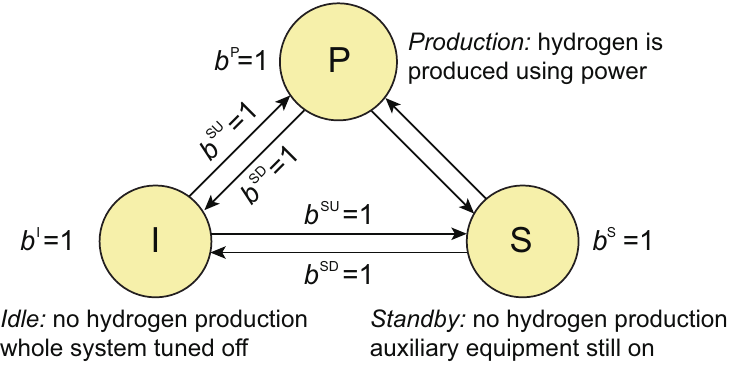}
  \caption{State transition diagram of the electrolyzer \cite{varela2021modeling,zheng2022optimal}.}
  \label{fig:states}
\end{figure}

The diagram of the state transition is shown in Fig. \ref{fig:states}, and a detailed description can be found in \cite{varela2021modeling} and \cite{zheng2022optimal}. Since the state transition model is not a contribution of this work, we will not elaborate it. We refer interested readers to the literature.

Note that several works assume a delay between \emph{Startup} and \emph{Production} \cite{li2021co}. 
Nevertheless, it is not a physical constraint. In practice, the electrolyzer is generally kept pressurized in \emph{Idle} to avoid the mechanical stress and energy cost of repressurization at \emph{Startup}. This allows for almost no delay between the power supply and hydrogen production. Indeed, the electrolyzers need to heat up to enable full-load operation. Because the temperature-related ramping limit is modeled in Section \ref{sec:temp}, we do not need a delay constraint in the state transition model.

\subsection{Hydrogen Production and Power Consumption}
\label{sec:prod}

The energy efficiency of an electrolyzer depends on the UI curve and Faraday efficiency \cite{ulleberg2003modeling}. The UI curve, also referred to as the \emph{polarization curve}, determines the overvoltage of the cell, which is usually approximated by:
\begin{align}
  U^{\text{cell}} & =  U^{\text{rev}} + (r_1 + r_2 T + r_3 P) I^{\text{cell}}  + s \log \left[ (t_1 + \frac{t_2}{T} + \frac{t_3}{T^2})  I^{\text{cell}} + 1 \right], \label{eq:ui}
\end{align}
\noindent
where $U^{\text{cell}}$ is the cell voltage, representing the electrical energy consumed by electrolysis; $U^{\text{rev}}$ is the reversible voltage, representing the energy converted into hydrogen; $T$ and $P$ are the temperature and pressure, respectively; and $s$, $r_1$, $r_2$, $r_3$, $t_1$, $t_2$, and $t_3$ are constant factors; see \cite{uchman2021varying} for details.

The electrolytic power $P_{i,j}^{\text{ele}}$ is determined by
\begin{align}
   P_{i,j}^{\text{ele}} = N^{\text{cell}} I^{\text{cell}} U^{\text{cell}}, \label{eq:pele}
\end{align}
\noindent
where $N^{\text{cell}}$ is the number of cells in the electrolyzer. The hydrogen production flow $\dot{n}^{\text{H}_2,\text{prod} }$ can be described by
\begin{align}
   \dot{n}^{\text{H}_2,\text{prod} } =  \frac{\eta^{\text{F}} N^{\text{cell}} I^{\text{cell}}}{2F}, \label{eq:nh2}
\end{align}
where $F=96485.3$ C/mol is the Faraday constant and $\eta^{\text{F}}$ is the Faraday efficiency of the electrolyzer, 
which is the percentage of electrons converted into product hydrogen \cite{ulleberg2003modeling}.

Suppose the pressure is maintained constant, a common practice to avoid fatigue \cite{qi2021pressure}. By combining (\ref{eq:ui})--(\ref{eq:pele}), the hydrogen production flow of an electrolyzer is a concave function of the electrolytic power and temperature \cite{li2021co,kopp2017energiepark}, denoted by
\begin{align}
   \dot{n}_{i,j}^{\text{H}_2,\text{prod} } = f(P_{i,j}^{\text{ele}},T_{i,j}). \label{eq:prod}
\end{align}

\begin{figure}[t]
  \centering
  \includegraphics[scale=1.15]{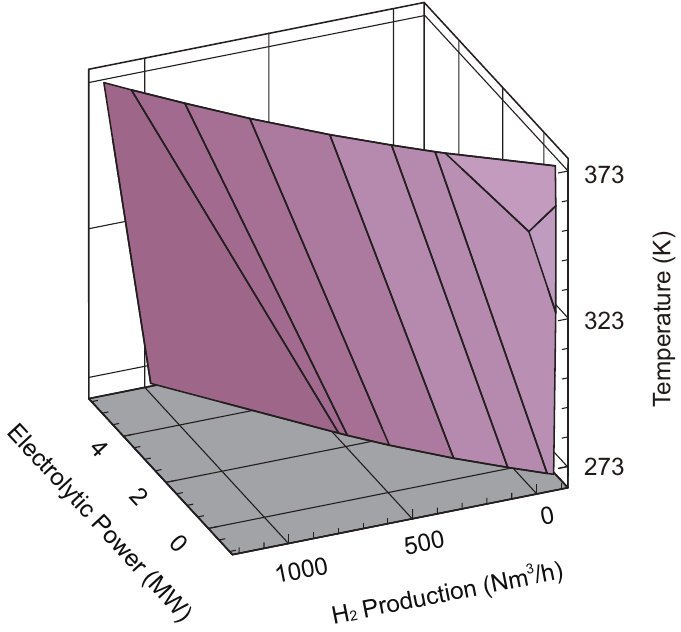}
  \caption{Approximate production function of a 1000 Nm$^3$/hr-rated electrolyzer.} 
  \label{fig:prod}
\end{figure}

To facilitate easy modeling of the plant scheduling problem without compromising the accuracy, the production function (\ref{eq:prod}) is approximated by a polyhedron using the double description (DD) algorithm \cite{jones2010polytopic} as shown in Fig. \ref{fig:prod} and then relaxed as a set of inequalities:
\begin{align}
  \dot{n}_{i,j}^{\text{H}_2,\text{prod} } \le \bm{A} P_{i,j}^{\text{ele}} + b_{i,j}^{\text{P}} T_{i,j} \bm{B}  +  b_{i,j}^{\text{P}} \bm{C}, \label{eq:prod}
\end{align}
\noindent
where $\bm{A}$, $\bm{B}$ and $\bm{C}$ are constant vectors. Hydrogen is only produced in the \emph{Production} state, expressed as
\begin{align}
  0 \le P_{i,j}^{\text{ele} } \le b_{i,j}^{\text{P}} \overline{P}_{i}^{\text{ele}}, \label{eq:flow}
\end{align}
\noindent
where $\overline{P}_{i}^{\text{ele}}$ is the dc power limit of the electrolyzer's rectifier.

Because the scheduling target always maximizes production, the operation point is forced onto the constraint surface of the production function, i.e., there is always at least an active constraint in (\ref{eq:prod}). In other words, we do not need binary variables to select the active subplane on the piecewise-linear function. 

\begin{remark} Many works adopt a linear production model for simplicity; see Table \ref{tab:literature}. However, linearization may lead to an error of 5\% \cite{zheng2022optimal}. To improve accuracy and fully exploit the flexibility of P2H production, this work considers the nonlinear efficiency.
\end{remark}


Moreover, to avoid sudden changes in the pressure, lye-gas separator liquid levels, and gas-liquid ratio in the electrolyzer that may cause excessive stress \cite{david2020dynamic}, the ramping of hydrogen production is limited as follows:
\begin{align}
  \underline{r}^{\text{H}_2,\text{prod} } \le \dot{n}_{i,j+1}^{\text{H}_2,\text{prod} }-\dot{n}_{i,j}^{\text{H}_2,\text{prod}} \le \overline{r}^{\text{H}_2,\text{prod}}, \label{eq:ramp}
\end{align}
\noindent
where $\overline{r}^{\text{H}_2,\text{prod}}$ and $\underline{r}^{\text{H}_2,\text{prod}}$ are the ramping bounds.

When fully heated, the ramping of an alkaline electrolyzer can reach as high as $\pm 20\%$ nominal load per second according to \cite{bertuccioli2014study} and our experiments shown in Appendix. Because the step length of scheduling is generally larger than 15 minutes, the ramping constraint (\ref{eq:ramp}) is not a hard limit. The temperature-related ramping limit is modeled in Section \ref{sec:temp}.

The power consumption of an electrolyzer is the sum of both electrolytic and auxiliary consumption:
\begin{align}
  P_{i,j} = P_{i,j}^{\text{ele}} + (b_{i,j}^{\text{P}} + b_{i,j}^{\text{S}}) P_{i,j}^{\text{aux}}, \label{eq:power}
\end{align}
\noindent where the auxiliary consumption $P_{i,j}^{\text{aux}}$ includes
\begin{align}
  P_{i,j}^{\text{aux}} = \frac{Q_{i,j}^{\text{heat}}}{\eta^{\text{heat}}}  +  \frac{Q_{i,j}^{\text{cool}}}{\eta^{\text{cool}}} + P^{\text{aux}}, \label{eq:bop}
\end{align}
\noindent
where $Q_{i,j}^{\text{heat}}$ and $Q_{i,j}^{\text{cool}}$ are heating and cooling heat (see Section \ref{sec:temp}); $\eta^{\text{heat}}$ and $\eta^{\text{cool}}$ are heating and cooling efficiencies; and $P^{\text{aux}}$ is the power of auxiliary equipment, including pumps, and control system, assumed to be constant here. 

\subsection{Temperature Dynamic Model and Constraints}
\label{sec:temp}

The temperature significantly impacts the electrolyzer efficiency and feasible load range \cite{grigoriev2020current,flamm2021electrolyzer,zheng2022optimal}. The impact comes from the electrolyte conductivity, activeness of the catalyst, and bubble coverage \cite{zheng2022optimal} and can be seen in (\ref{eq:prod}) and Fig. \ref{fig:prod}.
Unlike many previous works that omitted the temperature dynamics \cite{serna2017predictive,varela2021modeling,uchman2021varying}, this work includes it in the scheduling model to more precisely consider the temperature effect.

This work modified the famous first-order temperature model proposed by \cite{ulleberg2003modeling} 
to consider the state transitions in production scheduling. The block diagram is shown in Fig. \ref{fig:temp}(a), and the mathematical model is given below:
\begin{align}
  {C_i^\text{temp}} \frac{ T_{i,j+1}-  T_{i,j} }{h}  =  {Q_{i,j}^{\text{react}} - Q_{i,j}^{\text{diss}} - Q_{i,j}^{\text{cool}} + Q_{i,j}^{\text{heat}}}, \label{eq:temp}
\end{align}
\noindent
where $C^{\text{temp}}$ is the heat capacity of the electrolyzer; $T^{\text{am}}$ is the ambient temperature; $Q_{i,j}^{\text{react}}$ is the electrolytic reactional heat; $Q_{i,j}^{\text{diss}}$ and $Q_{i,j}^{\text{cool}}$ are the heat taken away by natural dissipation and active cooling, respectively; and $Q_{i,j}^{\text{heat}}$ is the auxiliary heating.

\begin{figure}[t]
  \centering
  \includegraphics[scale=1.15]{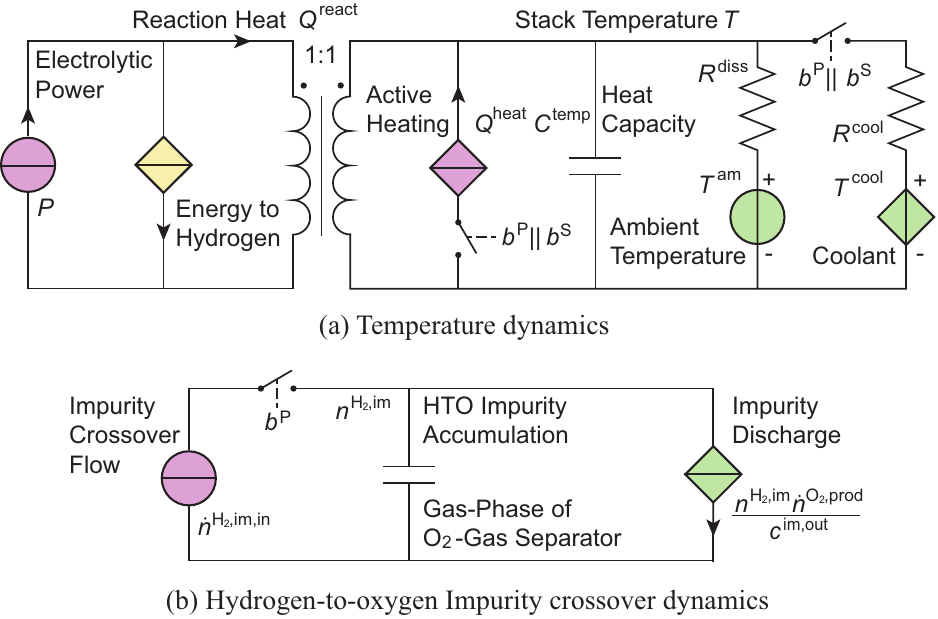}
  \caption{Equivalent circuits of the dynamic temperature and HTO impurity crossover model of the alkaline electrolyzer.}
  \label{fig:temp}
\end{figure}

The electrolytic heat $Q_{i,j}^{\text{react}}$ is determined as follows:
\begin{align}
  Q_{i,j}^{\text{react}} & = N^{\text{cell}} I_{i,j} \big(U_{i,j}^{\text{cell}} - U^{\text{th}} \big),  \label{eq:heatreactactual}
\end{align}
\noindent
where $U^{\text{th}}= 1.48$ V is the thermal neutral voltage and $I_{i,j}$ is the stack current, which satisfies (\ref{eq:nh2}). When the electrolytic voltage is higher than the thermal neutral voltage, which is a common situation in practice, excessive heat is produced. In fact, this is the main source of heating up the electrolyzer. The temperature affects $U_{i,j}^{\text{cell}}$ and therefore affects the electrolytic heat.

In the scheduling model, the electrolytic heat is approximated by a second-order function of the current and temperature:
\begin{align}
   Q_{i,j}^{\text{react}}  \approx N^{\text{cell}}  \big( a_0 I_{i,j} + a_1 I_{i,j} T_{i,j} + a_2 I^2_{i,j} -  U^{\text{th}} I_{i,j} \big), \label{eq:heatreact}
\end{align}
\noindent
where $a_0$, $a_1$ and $a_2$ are constant coefficients.

The natural heat dissipation $Q_{i,j}^{\text{diss}}$ to the environment follows
\begin{align}
  Q_{i,j}^{\text{diss}} = \frac{T_{i,j} - T^{\text{am}}}{R^{\text{diss}}}, \label{eq:heatdiss}
\end{align}
\noindent
where $R^{\text{diss}}$ is the thermal resistance of natural dissipation to the environment. The active cooling heat $Q_{i,j}^{\text{cool}}$ satisfies
\begin{align}
  0 \le P_{i,j}^{\text{cool}} \le  (b_{i,j}^{\text{P}} + b_{i,j}^{\text{S}}) \frac{T_{i,j} - T^{\text{cool}}}{R^{\text{cool}}}, \label{eq:cool}
\end{align}
\noindent
where $T^{\text{cool}}$ is the coolant temperature and $R^{\text{cool}}$ is the thermal resistance of cooling. The auxiliary heating $Q_{i,j}^{\text{heat}}$ is subject to
\begin{align}
  0\le Q_{i,j}^{\text{heat}} \le (b_{i,j}^{\text{P}} + b_{i,j}^{\text{S}}) \overline{Q}_{i,j}^{\text{heat}}, \label{eq:heat}
\end{align}
\noindent
where $\overline{Q}_{i,j}^{\text{heat}}$ is the upper limit. Note that auxiliary heating is available only in the \emph{Production} or \emph{Standby} state.

Although a higher temperature means a lower overvoltage and higher energy efficiency \cite{ulleberg2003modeling}, the stack temperature should stay below a limit. The temperature limit is twofold. First, to avoid damaging the diaphragm that separates the anode and cathode half-cells, the temperature should stay below a limit:
\begin{align}
  T_{i,j} \le \overline{T}\ (\text{usually}\ 85\ \text{to}\ 100\ ^\circ\text{C}). \label{eq:templimit}
\end{align}

Second, the electrolytic voltage should not exceed a safety margin (usually 2.1 V) to avoid damaging the microstructure of the electrode \cite{ulleberg2003modeling}. The constraint is approximated as follows:
\begin{align}
  U_{i,j}^{\text{cell}} \big( \approx a_0 + a_1 T_{i,j} + a_2 I_{i,j} \big) \le 2.1(V). \label{eq:voltlimit}
\end{align}

\begin{remark}
When the temperature is low, the cell voltage can be excessively high even at a relatively low load according to (\ref{eq:ui}). This contributes to the load ramping limit. In this work, because the temperature-related constraints are considered more precisely, we do not need a fixed ramping constraint.
\end{remark}

\subsection{Hydrogen-to-Oxygen Crossover Dynamics and Constraints}
\label{sec:hto}

\begin{figure}[t]
  \centering
  \includegraphics[scale=1.15]{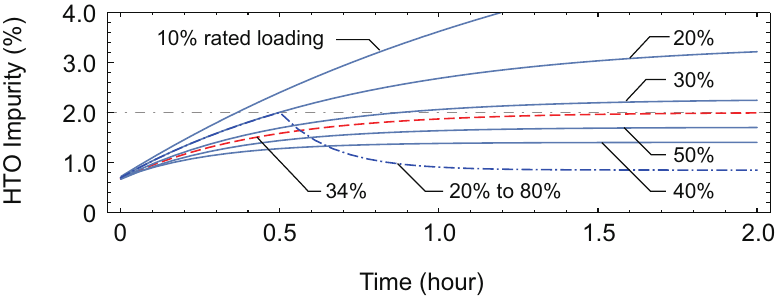}
  \caption{HTO crossover simulation of the alkaline electrolyzer with parameters given in Table \ref{tab:para} at different steady-state loading levels.}
  \label{fig:hto}
\end{figure}

Hydrogen-to-oxygen (HTO) impurity crossover may cause a flammable gas mixture. For safety, the electrolyzer shuts down when the hydrogen impurity in the oxygen product reaches 2\% in volume \cite{straka2021comprehensive}. A detailed description can be found in our previous work \cite{qi2021pressure}.
Due to the impurity accumulating faster at a low load, traditional scheduling models usually assume a lower limit of the feasible load level between 10\% and 40\% to avoid violating the limit \cite{varela2021modeling,uchman2021varying,he2021hydrogen}. However, this could conservatively limit the load range of the electrolyzer.

Because the accumulation of HTO impurity takes time, temporarily operating at a lower power level without violating the 2\% constraint is possible. This could extend the dynamic lower loading limit of the electrolyzer. Therefore, this work considers the dynamic impurity crossover limit.

The dynamic model of HTO impurity accumulation, proposed by our previous research \cite{qi2021pressure} and illustrated in Fig. \ref{fig:temp}(b), is as
\begin{align}
   \hspace{-2pt} n_{i,j+1}^{\text{H}_2,\text{im}} =  n_{i,j}^{\text{H}_2,\text{im}} + h \bigg( b_{i,j}^{\text{on}} \dot{n}^{\text{H}_2,\text{im},\text{in}} - \frac{ \dot{n}_{i,j}^{\text{O}_2,\text{prod}} {n}_{i,j}^{\text{H}_2,\text{im}}}{ c^{\text{im,out}}}  \bigg)  \label{eq:hto}
\end{align}
\noindent
where $\dot{n}^{\text{H}_2,\text{im},\text{in}}$ is the impurity crossover flow, assumed to be a constant because the system pressure is normally kept constant \cite{qi2021pressure}; $\dot{n}_{i,j}^{\text{O}_2,\text{prod}} = 0.5 \dot{n}_{i,j}^{\text{H}_2,\text{prod}}$ is the flow of product oxygen; and $c^{\text{im,out}}$ is a constant depicting the discharge rate of the impurity.

The HTO impurity constraint is finally expressed by limiting the percentage of hydrogen in the product oxygen as follows:
\begin{align}
    {n_{i,j}^{\text{H}_2,\text{im}} / n_{i,j}^{\text{O}_2,\text{prod}} } \le 2\%. \label{eq:conshto}
\end{align}

For ease of understanding, Fig. \ref{fig:hto} shows the accumulation of HTO impurities under different load levels with the electrolyzer parameters given in Table \ref{tab:para}. We can see that the critical steady-state lower limit is 34\%. Nevertheless, the load can stay lower for a short time without violating the 2\% safety limit. For example, the dot-dashed curve shows the impurity when the load is 20\% for the first 30 min and then raised to 80\%. Although it is initially loaded below the steady-state limit, the impurity does not exceed the limit. This allows for an extended load range compared to the traditional models with fixed load limits.




\begin{remark}
The dynamic temperature and impurity accumulation models introduced in Sections \ref{sec:temp} and \ref{sec:hto} are validated by experiments on an industry-rated electrolyzer. Detailed experimental settings and results are given in Appendix.
\end{remark}

\section{The Production Scheduling Problem and Decomposition-Based Solution Method}
\label{sec:scheduling}

\subsection{The P2H Plant Production Scheduling Model}
\label{sec:overall}

\subsubsection{Decision variables and scheduling horizon} The decision variables of the production scheduling are summarized in Table \ref{tab:constraints}. They determine the state transition and multiphysics dynamics of electrolyzers in the plant across a scheduling horizon $N$ with a step length of $h$. For clarity, the decision variables for the $i$th electrolyzer are denoted compactly as $\bm{x}_i$.

\subsubsection{Objective}
The objective is to maximize the profit of hydrogen production, expressed as the revenue from the sale of hydrogen minus the electricity cost and the startup cost that reflects the depreciation of the electrolyzers.

The objective can be separated for each electrolyzer. For the $i$th electrolyzer, the objective $\theta_i(\bm{x}_i)$ is minimizing:
\begin{align}
     \theta_i (\bm{x}_i) \triangleq \sum_{k=1}^{N} \left( -c^{\text{H}_{\text{2}}}  \dot{n}_{i,j}^{\text{H}_2,\text{prod} } + c^{\text{power}}_j P_{i,j} + c^{\text{SU}} b_{i,j}^{\text{SU}} \right), \label{eq:obj}
\end{align}
\noindent
where $c^{\text{H}_{\text{2}}}$ is the hydrogen selling price; $c^{\text{power}}_j$ is the electricity price at step $j$; and $c^{\text{SU}}$ is the startup cost. The overall objective is the sum of the $M$ electrolyzerwise objectives, as in (\ref{eq:objall}).

\subsubsection{Power supply constraints} Per the policies mentioned in the Introduction, we assume that 100\% renewable power is used. Hence, the P2H consumption cannot exceed the renewable energy source $P_k^{\text{{RES}}}$, denoted in vector form as
\begin{align}
  \bm{l}(\bm{p}) \triangleq \left[ \sum_{i=1}^{M}  P_{i,j} - P_j^{\text{{RES}}} \right]_{j=1,\ldots,N} \le \bm{0} , \label{eq:totalpower}
\end{align}
\noindent
where $P_j^{\text{{RES}}}$ is the available RES power supply at step $j$ and $\bm{p}$ is the vector of power consumption $P_{i,j}$ for all the $M$ electrolyzers at $N$ time steps. Note that $\bm{p}$ is the coupling vector of different electrolyzers and is separated as $\bm{p}=[\bm{p}_1^\mathrm{T},\bm{p}_2^\mathrm{T},\ldots,\bm{p}_M^\mathrm{T}]^\mathrm{T}$ with $\bm{p}_i=\bm{Q}_i \bm{x_i}$, where $\bm{Q}_i$ is a selection matrix.

\subsubsection{Operational feasibility constraints of the electrolyzers} The multiphysics dynamic constraints presented in Sections \ref{sec:prod} to \ref{sec:hto} are included in the scheduling model. Note that some of them have bilinear terms. We reformulate them to mixed-integer linear ones by the standard big-M method \cite{bemporad1999control}.

\begin{table}[tb]\footnotesize
  \renewcommand{\arraystretch}{1.3}
  \caption{Decision Variables and Multiphysics Feasibility Constraints of the Electrolyzer in the Scheduling Problem}\vspace{6pt}
  \label{tab:constraints}
  \centering
  \begin{tabular}{cccc}
   \hline \hline                         & \tabincell{c}{Decision variables $\bm{x}$}              & \tabincell{c}{Equality\\constr. $\bm{h}(\cdot)$\vspace{1.5pt}}  & \tabincell{c}{Inequality\\constr. $\bm{g}(\cdot)$}         \\ \hline
  State Transition     & $b^{\text{P}}, b^{\text{S}},b^{\text{I}},b^{\text{SU}}$                  & (\ref{eq:states})                                   & (\ref{eq:startup}), (\ref{eq:gap})               \\
  Production           & $\dot{n}^{\text{H}_2,\text{prod}}, I, P, P^{\text{ele}}, P^{\text{aux}}$    & (\ref{eq:nh2}), (\ref{eq:power}), (\ref{eq:bop})    & (\ref{eq:prod})--(\ref{eq:ramp})                           \\
  Temperature          & $T,Q^{\text{react}},Q^{\text{diss}},Q^{\text{cool}},Q^{\text{heat}}$     & (\ref{eq:temp}), (\ref{eq:heatreact})--(\ref{eq:heatdiss})                & (\ref{eq:cool})--(\ref{eq:voltlimit})                  \\
  HTO Impurity         & $n^{\text{H}_2,\text{im}}$                                               & (\ref{eq:hto})                                      & (\ref{eq:conshto})                      \\
  \hline \hline
  \end{tabular}
\end{table}

We denote the feasible operation region of the $i$th electrolyzer subject to the multiphysics constraints as $\Omega_i$, written as
\begin{align}
  \bm{x}_i \in \Omega_i \triangleq \big\{ \bm{x}_i :
      \bm{g}_i(\bm{x}_i)  \le \bm{0} ,
      \bm{h}_i(\bm{x}_i)  = \bm{0} \big\},  \label{eq:conselyzers}
\end{align}
\noindent
where the entries of $\bm{h}_i(\cdot)$ and $\bm{g}_i(\cdot)$ are summarized in Table \ref{tab:constraints}. For different electrolyzers, the constraints are mutually independent.

\subsubsection{Overall plant production scheduling problem and its scalability issue}
Summarizing the above, the overall P2H plant production scheduling problem is formulated as follows:

{\vspace{3pt}\noindent\bf Plant scheduling problem (PSP):}
\begin{align}
  \min_{\bm{x}_1,\ldots,\bm{x}_M}  \hspace{10pt} &  \theta_1 (\bm{x}_1) + \theta_2 (\bm{x}_2) + \ldots + \theta_{M} (\bm{x}_{M}) \label{eq:objall}  \\
  \text{s.t.} \hspace{20pt} &   \bm{x}_i \in \Omega_i,\ \forall i=1,\ldots,M, \label{eq:consall} \\
  &
  \bm{l}(\bm{p}) \le \bm{0},  \label{eq:supplyall}  
\end{align}
\noindent
where the components in the objective (\ref{eq:objall}) are given in (\ref{eq:obj}); the electrolyzer operation feasibility constraints (\ref{eq:consall}) are defined in (\ref{eq:conselyzers}); and the power supply constraint (\ref{eq:supplyall}) is given as (\ref{eq:totalpower}).

The PSP (\ref{eq:objall})--(\ref{eq:supplyall}) is a large-scale and nonconvex MILP. For instance, given the step length $h=15$ min and horizon $N=96$, i.e., a day, an electrolyzer has $2,208$ real decision variables and $869$ binary decision variables, i.e., $\mathrm{dim}(\bm{x}_i)=3,077$. For a plant comprised of $20$ electrolyzers, we have $61,540$ decision variables, and its nonconvexity is stronger than that of the traditional scheduling methods due to the multiphysics dynamics. It does not converge for days if we solve it directly using off-the-shelf solvers such as \emph{Gurobi}. Therefore, an efficient solution method is needed to enable the scheduling of an industry-scale P2H plant.

\subsection{Parallelizable Decomposition-Based Solution Method}
\label{sec:decomposition}

The PSP (\ref{eq:objall})--(\ref{eq:supplyall}) has a natural separable structure, which is helpful for developing an efficient solution method of the plant production scheduling problem. The objective (\ref{eq:objall}) and the electrolyzer operation feasibility constraints (\ref{eq:consall}) can be separated with respect to each electrolyzer. The only coupling is the power supply (\ref{eq:supplyall}). Thus, the large-scale nonconvex problem can be decomposed into $M$ MILP subproblems, each representing one electrolyzer.

Nevertheless, because the subproblems are nonconvex, classical decomposition-based methods such as the 
ADMM fails to converge. To address the nonconvexity while leveraging the separable structure, we adopt a recent algorithm called SDM-GS-ALM \cite{boland2019parallelizable}, which is a combination of the simplicial decomposition method (SDM), block Gauss-Seidel method (GS), and augmented Lagrangian multiplier method (ALM).

Specifically, the augmented Lagrangian (AL) $L_{\rho,i}$ for subproblem $i(=1,\ldots,M)$ with respect to electrolyzer $i$ is defined:
\begin{align}
  L_{\rho,i} (\bm{x}_i,\bm{p}_i,\bm{\omega}_i) & \triangleq \theta_i(\bm{x}_i)  
  + \bm{\omega}_i^\mathrm{T} (\bm{Q}_i \bm{x}_i - \bm{p}_i) + \frac{\rho}{2} \Vert \bm{Q}_i \bm{x}_i - \bm{p}_i  \Vert_2^2  \nonumber 
\end{align}
\noindent
where $\bm{\omega}_i$ is the Lagrange multiplier and $\rho$ is the penalty.

Instead of the nonconvex feasibility region $\Omega_i$, the subproblem is first solved in a polyhedron $D_i \subseteq \mathrm{conv}(\Omega_i)$. Thus, the subproblems are relaxed as convex. Initially, $D_i$ is set as a single feasible point, i.e., $D_i = \{\bm{x}_i^0\}$. In the $k$th round iteration of SDM-GS-ALM, we solve the relaxed subproblems (RSPs) and a coupling problem (CP) alternatively to update $\bm{x}_i$ and the coupling variable $\bm{p}$, known as the block GS method.

{\vspace{2pt}\noindent\bf Relaxed Subproblem $i$ (RSP$_i$):}
\begin{align}
  \tilde{\bm{x}}_i = \arg\min_{\bm{x}_i} \big\{  L_{\rho,i}&  (\bm{x}_i,\tilde{\bm{p}}_i,\bm{\omega}_i^{k}) : \bm{x}_i \in D_i \big\}, \label{eq:rsp}
\end{align}

{\vspace{2pt}\noindent\bf Coupling Problem (CP):}
\begin{align}
  \tilde{\bm{p}} = \arg\min_{\bm{p}} \Big\{ \sum\nolimits_{i=1}^{M} L_{\rho,i}&  (\tilde{\bm{x}}_i,\bm{p}_i,\bm{\omega}_i^{k}) : \bm{l}(\bm{p}) \le \bm{0} \Big\}. \label{eq:cp}
\end{align}

After $t_{\text{max}}$ rounds of alternations of (\ref{eq:rsp}) and (\ref{eq:cp}), we update the decision variables as $\bm{x}_i^k \leftarrow \tilde{\bm{x}}_i$ and $\bm{p}^k \leftarrow \tilde{\bm{p}}$. Then, we find a new point $\hat{\bm{x}}_i$ in the feasibility region $\Omega_i$ by solving:

{\vspace{2pt}\noindent\bf Convex Hull Update Problem $i$ (CHUP$_i$):}
\begin{align}
  \hat{\bm{x}}_i \in \arg\min_{\bm{x}_i} \big\{ (\bm{x}_i-\bm{x}^k_i)^{\mathrm{T}}  \nabla_{\bm{x}_i}  L_{\rho,i} : \bm{x}_i \in \Omega_i \big\},  \label{eq:ppp}
\end{align}
\noindent
and the polyhedron $D_i$ is updated by the convex hull operation:
\begin{align}
  D_i \leftarrow \mathrm{conv}(D_i,\hat{\bm{x}}_i). \label{eq:updateconvexhull}
\end{align}

The term SDM appears in the name of the SDM-GS-ALM algorithm because the feasibility regions $D_i$ of the subproblems are polyhedrons, i.e., a simplicial set. The nonconvexity of the subproblem is thus far avoided. 

\begin{algorithm}[tb]
  \caption{SDP-GS-ALM for P2H Production Scheduling}
  \label{alg:alm}
  \begin{algorithmic}[1]
    \REQUIRE $\gamma \in (0,1)$, $t_{\text{max}}$, $\epsilon>0$, $\rho>0$
    \STATE initialize $k \leftarrow 0$, $\check\phi^{0} \leftarrow \infty $, $\bm{x}^0_i \in \Omega_i$, $D_i\leftarrow\{\bm{x}^0_i\}$
    \REPEAT
    \STATE $k \rightarrow k+1$
    \STATE $\bm{\omega}_i^k \leftarrow \bm{\omega}_i^{k-1}$, $\tilde{\bm{x}}_i \leftarrow \bm{x}_i^{k-1}$, $\tilde{\bm{p}} \leftarrow \bm{p}^{k-1}$, $\check\phi^{k} \leftarrow \check\phi^{k-1}$
    \FOR{$t=1,2,\ldots,t_{\text{max}}$}
      \STATE update $\tilde{\bm{x}}_i$ by the RSP$_i$ (\ref{eq:rsp}) for $i=1,2,\ldots,M$
      \STATE update $\tilde{\bm{p}}$ the CP (\ref{eq:cp})
    \ENDFOR
    \STATE find $\hat{\bm{x}}_i$ by the CHUP$_{i}$ (\ref{eq:ppp}), for $i=1,2,\ldots,M$
    \STATE $D_i^k \leftarrow \mathrm{conv}(D_i^{k-1},\hat{\bm{x}}_i)$, for $i=1,2,\ldots,M$
    \STATE examine the serious step condition by (\ref{eq:ssc1}) and (\ref{eq:ssc})
    \IF{$\eta_k > \gamma$}
      \STATE $\bm{\omega}^{k}_i \leftarrow \bm{\omega}^{k-1}_i + \rho (\bm{Q}_i \bm{x}^k_i - \bm{p}^k_i)$, for $i=1,2,\ldots,M$
      \STATE $\check{\phi}^k \leftarrow \sum_{i=1}^{M} \tilde\phi_i$
    \ENDIF
    \UNTIL{$\sum_{i=1}^{M} \big( L_{\rho,i}(\bm{x}^{k}_i,\bm{p}^{k}_i,\bm{\omega}^{k-1}_i) + \frac{\rho}{2} \Vert \bm{Q}_i \bm{x}^k_i - \bm{p}^k_i \Vert^2_2 \big)/ \hat{\phi}^{k}<\epsilon$}
  \end{algorithmic}
\end{algorithm}

To allow for a larger penalty $\rho$ and fast convergence while mitigating the destabilizing effect, the serious step condition (SSC) \cite{boland2019parallelizable} is examined at each iteration of $k$, as follows:
\begin{align}
  \tilde\phi_i &=  L_{\rho,i} (\bm{x}^{k}_i,\bm{p}^{k}_i,\bm{\omega}^{k}_i) + \frac{\rho}{2} \Vert \bm{Q}_i \bm{x}^k_i - \bm{p}^k_i \Vert^2_2  - \Gamma_i, \label{eq:ssc1}\\
  \eta^k & = \frac{\sum_{i=1}^{M} \tilde\phi_i  - \check\phi^{k}}{ \sum_{i=1}^{M} \big( L_{\rho,i}(\bm{x}^{k}_i,\bm{p}^{k}_i,\bm{\omega}^{k}_i) + \frac{\rho}{2} \Vert \bm{Q}_i \bm{x}^k_i - \bm{p}^k_i \Vert^2_2 \big) - \check\phi^{k}},  \label{eq:ssc}
\end{align}
\noindent
where $\Gamma_i =  \nabla_{\bm{x}_i}  L_{\rho,i} (\tilde{\bm{x}}_i,\tilde{\bm{p}}_i,\bm{\omega}^{k}_i)^{\mathrm{T}} (\hat{\bm{x}}_i-\tilde{\bm{x}}_i)$.
If $\eta_k > \gamma$, the multiplier $\bm{\omega}_i$ is updated by
\begin{align}
  \bm{\omega}^k_i \leftarrow \bm{\omega}^{k}_i + \rho (\bm{Q}_i \bm{x}^k_i - \bm{p}^k_i). \label{eq:updatew}
\end{align}

In each iteration of $k$, with possibly a new vertex added, the polyhedron $D_i$ expands and is ensured to include the optimal solution of the PSP \cite{boland2019parallelizable}. The overall procedure of the solution method is summarized as Algorithm \ref{alg:alm}.

\begin{remark}
  Among the decomposed subproblems (RSP, CP, and CHUP), only the CHUP (\ref{eq:ppp}) is nonconvex. It is related to only one electrolyzer. The other subproblems are convex and easy to solve. Therefore, the large-scale scheduling problem is replaced by small-scale problems via the decomposition algorithm, and the scalability issue mentioned above is solved. Numerical examples can be seen in Section
\end{remark}

The SDM-GS-ALM method is also used in other engineering problems, such as power system unit commitment (UC) \cite{chen2019decentralized} and integrated energy system (IES) scheduling \cite{madadi2019decentralized}, to which interested readers are referred to.

\begin{table}[tb]\footnotesize
  \renewcommand{\arraystretch}{1.3}
  \caption{Parameters of the Electrolyzer Used in the Case Study}\vspace{6pt}
  \label{tab:para}
  \centering
  \begin{tabular}{ll}
  \hline \hline
  Parameter    \hspace{86pt}                                                    & Value \hspace{60pt}                 \\  \hline
  Rated hydrogen production                                                     & $1,000$ Nm$^3$/h$^2$ \\
  Maximal power $\overline{P}^{\text{ele}}$                                     & $6$  MW \\
  Number of cells $N^{\text{cell}}$                                             & $260$   \\
  Ramping limits $\overline{r}^{\text{H}_2,\text{prod}} / \underline{r}^{\text{H}_2,\text{prod} }$    & $+1,600 / -4,800$ Nm$^3$/h$^2$         \\
  Production function $f(\cdot)$                                                & See Fig. \ref{fig:prod}                   \\ \hline
  Coolant temperature $T^{\text{cool}}$                                         & $278$ K ($5\ ^\circ$C) \\
  Temperature limit $\overline{T}$                                              & $368$ K ($95\ ^\circ$C) \\
  Heat capacity $C^{\text{temp}}$                                               & $1.163\times10^8$ J/K    \\
  Dissipation resistance $R^{\text{diss}}$                                      & $1.2\times10^{-4}$ K/W    \\
  Active cooling resistance $R^{\text{cool}}$                                   & $2\times10^{-5}$ K/W    \\\hline
  Impurity crossover flow $\dot{n}^{\text{H}_2,\text{im},\text{in}}$           & $0.003182$ mol/s  \\
  Impurity discharge constant $c^{\text{im,out}}$                               & $5.68\times10^5$ mol$^{-1}$ \\
  \hline \hline
  \end{tabular}
\end{table}


\section{Case Studies}
\label{sec:case}

\subsection{Case Settings}
\label{sec:setting}


The proposed P2H scheduling method is tested with two plant settings. The first plant comprises $4$ electrolyzers, each rated $5$ MW ($1,000$ Nm$^3$/h), and is used to demonstrate the impact of the multiphysics effects. The second has $22$ electrolyzers, with a total rating of $110$ MW ($22,000$ Nm$^3$/h), and is used to exhibit the proposed method's ability to deal with large-scale problems.

The horizon of production scheduling is set as one day, i.e., $N^h = 96$, with a step length $h$ of $15$ min. Without loss of generality, the hydrogen price and electricity price are set as constant, i.e., $0.38$ \$/Nm$^3$ and $34.7$ \$/MWh, respectively. The parameters of the electrolyzer are given in Table \ref{tab:para}. The efficiencies of the electrolyzers in a plant are perturbed within $5\%$ to simulate different degradation statuses. The startup cost is set as $280$ \$. The simulation platform is \emph{Wolfram Mathematica 12.3}, and the optimization solver employed is \emph{Gurobi 9.5.0}.

\subsection{Base-Case Scheduling Result of a 4-Electrolyzer Plant}
\label{sec:result}

\begin{figure}[t]
  \centering
  \includegraphics[scale=1.15]{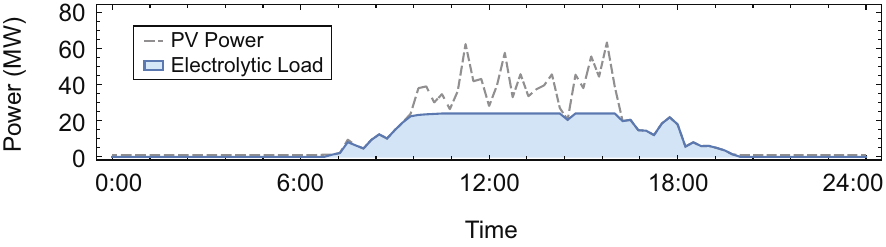}
  \caption{Power supply and load of the 4-electrolyzer plant under the proposed production scheduling method.}
  \label{fig:power}
\end{figure}

\begin{figure}[t]
  \centering
  \includegraphics[scale=1.15]{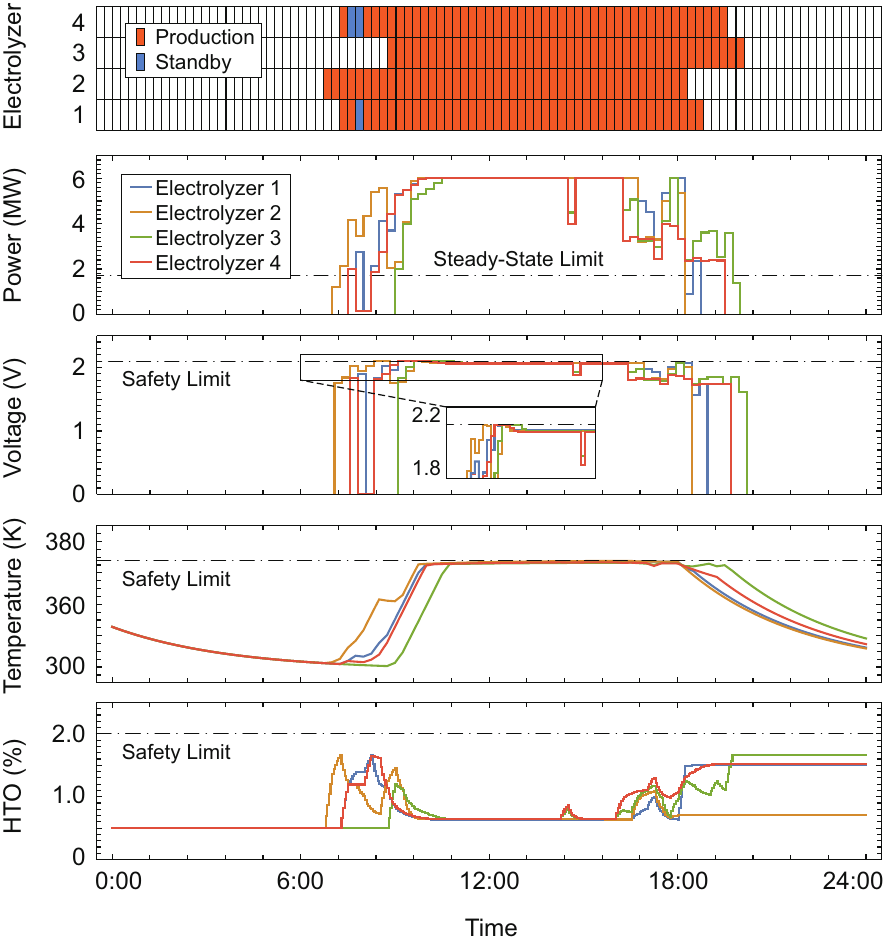}
  \caption{Electrolyzer state transition, power, cell voltage, temperature, and HTO impurity of the 4-electrolyzer plant under the proposed scheduling method.}
  \label{fig:multi}
\end{figure}

We assume that the P2H plant is directly connected to photovoltaic power. The power supply profile is based on the data of a PV plant in Sichuan Province, China, as shown in Fig. \ref{fig:power} \cite{qiu2020stochastic}.
The proposed scheduling method calculates the state transition and power of the four electrolyzers. Then, by time-domain simulation of the detailed dynamic thermal and mass transfer models, the electrolyzer cell voltage, temperature, and HTO impurity are exhibited in Fig. \ref{fig:multi}.

As observed, the electrolyzers start up following sunrise and shut down after sunset. Due to the temperature-related cell voltage limit (\ref{eq:voltlimit}), the load of each electrolyzer increases gradually as the temperature increases. After fully heating up, the voltage drops slightly below the limit even if the electrolyzers are fully loaded. Most of the time, the electrolyzers operate at the upper limit of temperature to maximize production.

At startup, the HTO impurity first increases due to a temporary low load. At approximately 7:30, due to a temporary drop in PV output, two electrolyzers switch to \emph{Standby} to avoid a long-term low load and exceeding the HTO impurity limit. At other times, the HTO impurity is relatively low due to the high load. Although the steady-state load limit is 34\%, 
the proposed scheduling method enables the electrolyzers to operate at a lower load temporarily. This allows more renewable energy to be utilized compared to the traditional scheduling method \cite{varela2021modeling}.


\begin{table}[tb]\footnotesize
  \renewcommand{\arraystretch}{1.35}
  \caption{Hydrogen Output and Profit of Different Scheduling Methods Considering Different Multiphysics Constraints}\vspace{6pt}
  \label{tab:comparison}
  \centering
  \begin{tabular}{ccccc}
  \hline \hline
  Method                                                                    & \tabincell{c}{Hydrogen\\Output (Nm$^3$)}                          & \tabincell{c}{Electricity\\Cost (\$)}                         & \tabincell{c}{Startup\\Cost (\$)}      & Profit (\$)                                                    \\ \hline
  \tabincell{c}{Traditional \cite{varela2021modeling}\\(w/o multiphysics)} & $39760.9$                                                         & $7867.0$                                                      & $1120.0$                               & $6122.2$                                                       \\ \vspace{-9pt}  \\
  \tabincell{c}{Traditional\\+temp. dynamics}                               & \tabincell{c}{$39822.9$\vspace{1.5pt}\\($+0.15\%$)}               & \tabincell{c}{$\bm{7818.7}$\vspace{1.5pt}\\($\bm{-0.61\%}$)}  & \tabincell{c}{$1120.0$\vspace{1.5pt}}  & \tabincell{c}{$6194.0$\vspace{1.5pt}\\($+1.17\%$)}             \\ \vspace{-7.5pt}  \\
  \tabincell{c}{Traditional\\+HTO dynamics}                                 & \tabincell{c}{$39912.7$\vspace{1.5pt}\\($+0.38\%$)}               & \tabincell{c}{$7901.6$\vspace{1.5pt}\\($+0.04\%$)}            & \tabincell{c}{$1120.0$\vspace{1.5pt}}  & \tabincell{c}{$6145.2$\vspace{1.5pt}\\($+0.38\%$)}             \\ \vspace{-7.5pt}  \\
  \tabincell{c}{\bf Proposed (w/ temp.\\\bf\&HTO dynamics)}                 & \tabincell{c}{$\bm{40024.0}$\vspace{1.5pt}\\($\bm{+0.825\%}$)}    & \tabincell{c}{$7854.7$\vspace{1.5pt}\\($-0.16\%$)}            & \tabincell{c}{$1120.0$\vspace{1.5pt}}  & \tabincell{c}{$\bm{6234.4}$\vspace{1.5pt}\\($\bm{+1.83\%}$)}   \\ \vspace{-11pt}    \\
  \hline \hline
  \end{tabular}
\end{table}

\begin{figure}[t]
  \centering
  \includegraphics[scale=1.15]{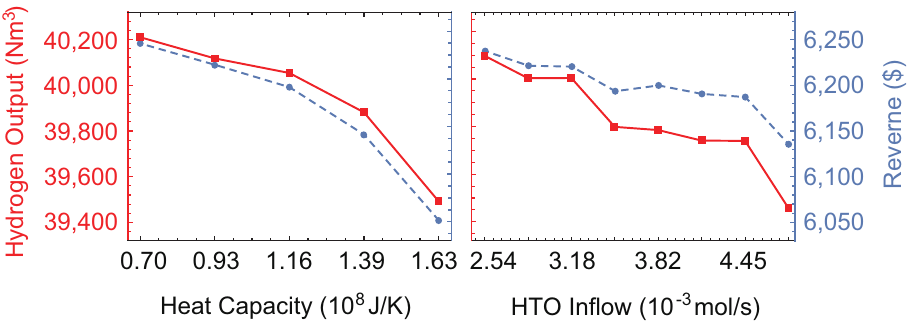}
  \caption{Impact of the electrolyzer heat capacity and HTO impurity inflow rate on the hydrogen output and revenue of daily hydrogen production.}
  \label{fig:physics}
\end{figure}

\subsection{Significance of Considering Multiphysics Effects}
\label{sec:compare}

To exhibit the significance of considering the dynamic temperature and HTO impurity effects, we quantitatively compare the proposed multiphysics-aware scheduling approach to the traditional method \cite{varela2021modeling} under the same plant setting and power supply. For completeness, we also separately add temperature and HTO dynamics to the traditional scheduling model. The results of hydrogen output, electricity cost, startup cost, and total profit are presented in Table \ref{tab:comparison}.

We can see that hydrogen production increases with less electricity consumed after considering the temperature effect. In other words, the energy conversion efficiency is improved. Considering the HTO impurity effect, hydrogen production increases with more electricity consumed. This indicates that extra renewable energy, especially when the power supply is low and the steady-state limit does not allow production by the electrolyzer, can be utilized compared to the traditional method.

Comparatively, as the proposed method considers both the temperature and HTO impurity dynamics, the hydrogen output increases by $0.825\%$, and the total profit increases by $1.83\%$. This is larger than the sum of the improvements of $1.17\%$ and $0.38\%$ by either independently considering the temperature or HTO impurity effects. An extra benefit of coordinating the multiphysics effects in the optimal scheduling of hydrogen production is therefore observed.


\subsection{Sensitivity Analysis of the Multiphysics Parameters}

To show the impact of the multiphysics parameters of the electrolyzer on hydrogen production, based on the settings of Section \ref{sec:result}, the results of different heat capacities and HTO impurity inflow rates are presented in Fig. \ref{fig:physics}. As observed, the electrolyzers heat up faster with the decrease in the heat capacity. This allows for faster ramping and higher efficiency, resulting in higher hydrogen output. However, if the energy source is highly fluctuating, a larger heat capacity better retains the temperature and thus is conversely more adaptable to a flexible operation.

As for HTO impurity accumulation, with the increase in the impurity inflow rate, the hydrogen output and profit drop. This phenomenon is attributed to the lower-power operation of the electrolyzers being further limited. The above result presents a quantitative analysis of the motivation for reducing HTO crossover, such as improving the diaphragm's performance.


\begin{figure}[t]
  \centering
  \includegraphics[scale=1.15]{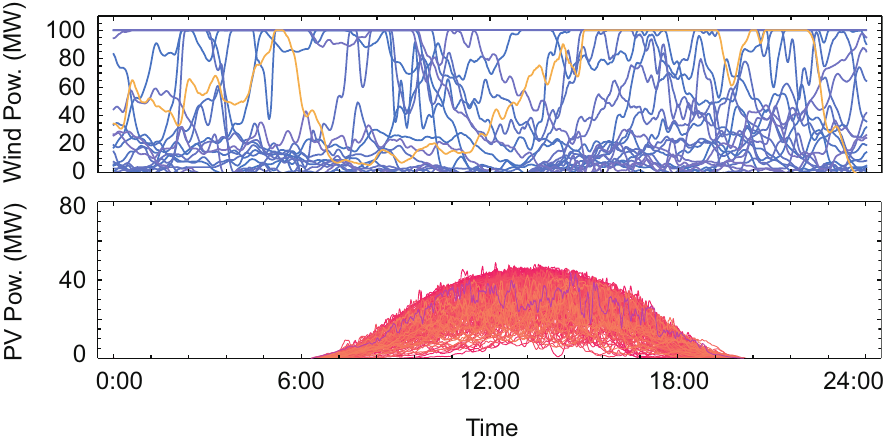}
  \caption{Wind and PV power scenarios to test the P2H scheduling method.}
  \label{fig:scen}
  \vspace{12pt}
\end{figure}

\begin{figure}[t]
  \centering
  \includegraphics[scale=1.15]{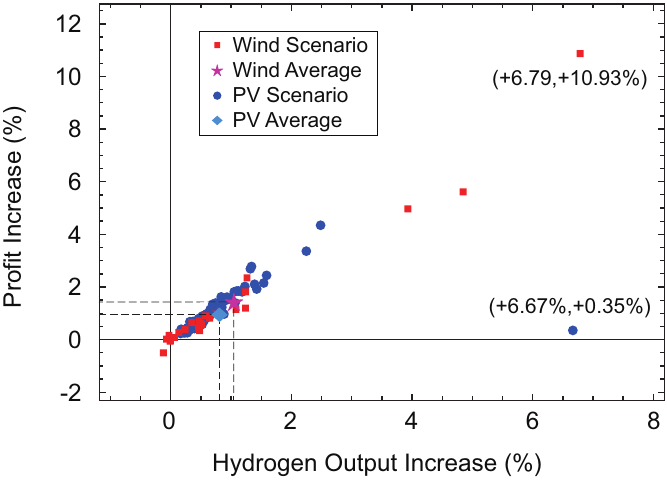}
  \caption{Improvements of the proposed multiphysics dynamics-aware scheduling method compared to the traditional scheduling method in terms of hydrogen output and total revenue under various wind and PV generation scenarios.}
  \label{fig:delta}
\end{figure}

\subsection{Comparison Under Various Wind and PV Scenarios}
\label{sec:scen}

The proposed and the traditional P2H scheduling methods are compared under various wind and PV power supply scenarios, as shown in Fig. \ref{fig:scen}. The $25$ wind scenarios are based on an offshore wind farm collected by Ris\o{} \cite{lin2012assessment}, and the $100$ PV scenarios are based on a plant in Sichuan Province, China \cite{qiu2020stochastic}.

The simulation results, in terms of the increases in hydrogen output and profit of the proposed multiphysics-aware method compared to the traditional method, are plotted in Fig. \ref{fig:delta}. We can observe that under both wind and PV power supply scenarios, the proposed method achieves higher hydrogen output and profit. On average, with wind power, the proposed method leads to a $1.021\%$ increase in hydrogen output and a $1.438\%$ increase in profit. For the PV scenarios, the average increases are $0.807\%$ and $0.982\%$. Considering the large investment for P2H projects, these improvements are significant.

\begin{figure}[t]
  \centering
  \includegraphics[scale=1.15]{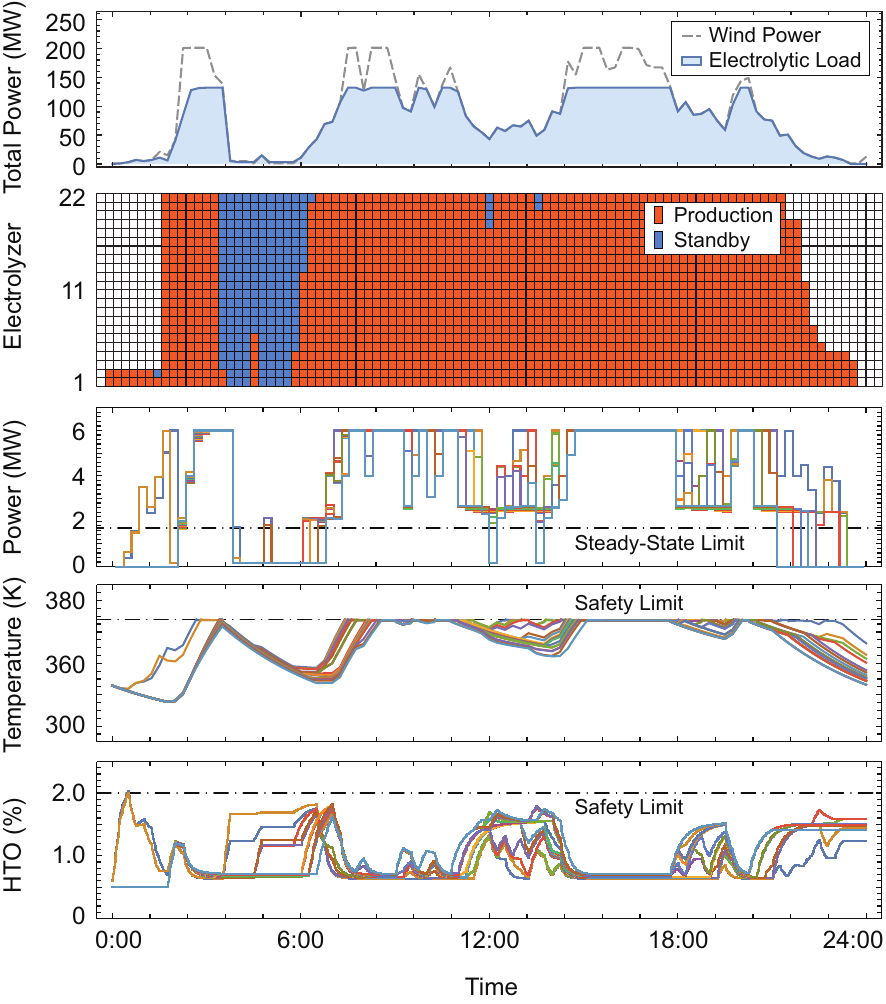}
  \caption{Power supply, electrolyzer state transition, power, temperature, and HTO impurity of the 22-electrolyzer plant obtained by the proposed method.}
  \label{fig:multilarge}
\end{figure}

\subsection{Scheduling of a Large-Scale Plant with 22 Electrolyzers}
\label{sec:large}

To show the proposed method's ability to schedule the production of a large-scale P2H plant, we test it using a plant setting of 22 electrolyzers. The scale is based on a real-life P2H plant under construction in Inner Mongolia, China, and we assume it is directly connected to wind power.

Using the proposed method in Section \ref{sec:decomposition}, the tolerance is set as $\epsilon = 0.1\%$, $\mathrm{mipgap}$ for solving the subproblems is $ 10^{-6}$, and the computation time is $713.82$ s. The supply and total electrolytic power and the electrolyzer state transition, load, temperature, and HTO impurity are plotted in Fig. \ref{fig:multilarge}.
Similar to the $4$-electrolyzer plant case, the temperature and HTO impurity constraints are pushed to the preset limits to maximize the flexibility of the electrolyzer. The total hydrogen production is $305,709.3$ Nm$^3$, and the profit is $47,809.9$ \$. Compared to the traditional scheduling method \cite{varela2021modeling}, the hydrogen production is increased by $7.74\%$, and the profit is increased by $8.72\%$. Due to space limitations, a detailed comparison is not given here.

\subsection{Analysis of the Scalability of the Solution Methods}

\begin{figure}[t]
  \centering
  \includegraphics[scale=1.15]{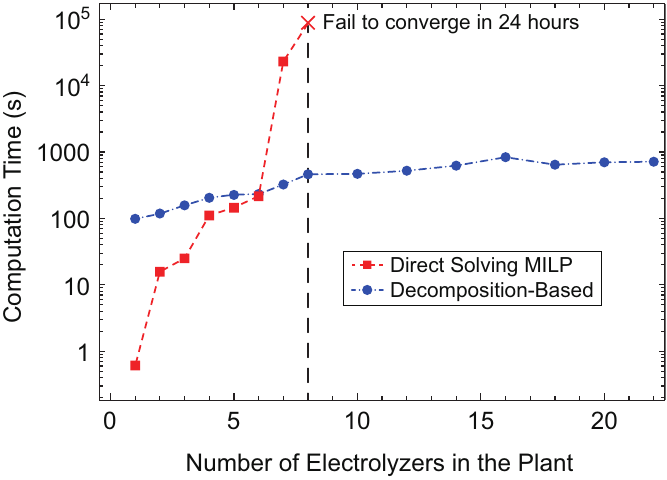}
  \caption{Computation time of directly solving the plant scheduling problem via the proposed decomposition-based method.}
  \label{fig:time}
\end{figure}

Finally, we compare the computation time of directly solving the plant scheduling problem or using SDM-GS-ALM proposed in Section \ref{sec:decomposition}. For SDM-GS-ALM, the computation parameters are the same as in Section \ref{sec:large}. When directly solving the PSP using \emph{Gurobi}, $\mathrm{mipgap}$ is set as $0.1\%$. The computation time with different plant sizes is shown in Fig. \ref{fig:time}. If we solve directly, the computation time increases exponentially. For $8$ electrolyzers, it fails to converge in $24$ hours. Nevertheless, the proposed decomposition-based solution method has linear complexity. For $8$ electrolyzers, it takes $461.68$ s, and for $22$ electrolyzers, it takes only $713.82$ s. This makes it feasible to schedule a large-scale plant with tens of electrolyzers, considering the complicated multiphysics dynamics.

\section{Conclusions}

This paper first incorporates experimentally validated multiphysics dynamic models of AEL into the production scheduling of industry-scale P2H plants.
A decomposition-based solution method is then proposed to offset the increased complexity. The case studies show that by considering the multiphysics effects, the loading flexibility of the P2H plant is significantly improved, leading to an average profit improvement of $1.438\%$ or $0.982\%$ when the plant is directly coupled with wind or solar energy.

The current work relies on an accurate wind and solar energy forecast and does not consider forecast errors. Considering the uncertainty of renewable power, a receding-horizon dispatch framework or a stochastic optimization method to alleviate the impact of uncertainty is needed in future studies.

In addition, industrial P2H production can be coupled with chemical plants, e.g., ammonia or methanol synthesis, as shown in Fig. \ref{fig:cluster}. Therefore, further consideration of the multiphysics dynamic limits of chemical plants in a joint power-to-chemicals (P2X) framework is one of the future works.



\section*{Appendix: Experimental Validation of the Multiphysics Model}
\label{sec:app1}

\begin{figure}[t]
  \centering
  \includegraphics[width=3.4in]{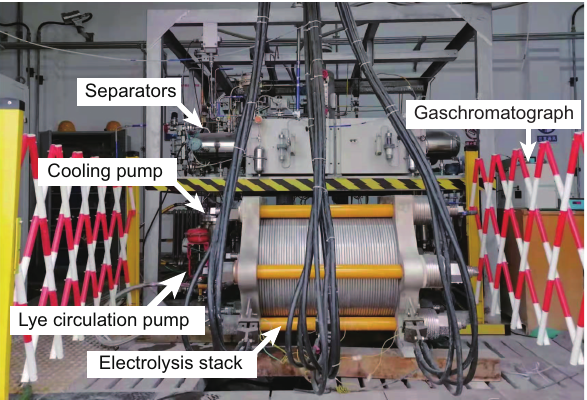}
  \caption{The CNDQ5/3.2 alkaline electrolyzer used for model validation.}
  \label{fig:photo}
\end{figure}

The dynamic models of the temperature and HTO impurity effects presented in Sections \ref{sec:temp} and \ref{sec:hto} are verified by an experiment on a CNDQ5/3.2 alkaline electrolyzer, which is manufactured by the Purification Equipment Research Institute of China Shipbuilding Industry Corporation (CSIC), as shown in Fig. \ref{fig:photo}. Its power rating is $25$ kW, and the rated hydrogen flow is $5$ Nm$^3$/h. The rescaled PJM RegD regulation signal \cite{pjmsginal} on Dec. 1, 2019, is used as the power command of the electrolytic power.

\begin{figure}[t]
  \centering
  \includegraphics[scale=1.15]{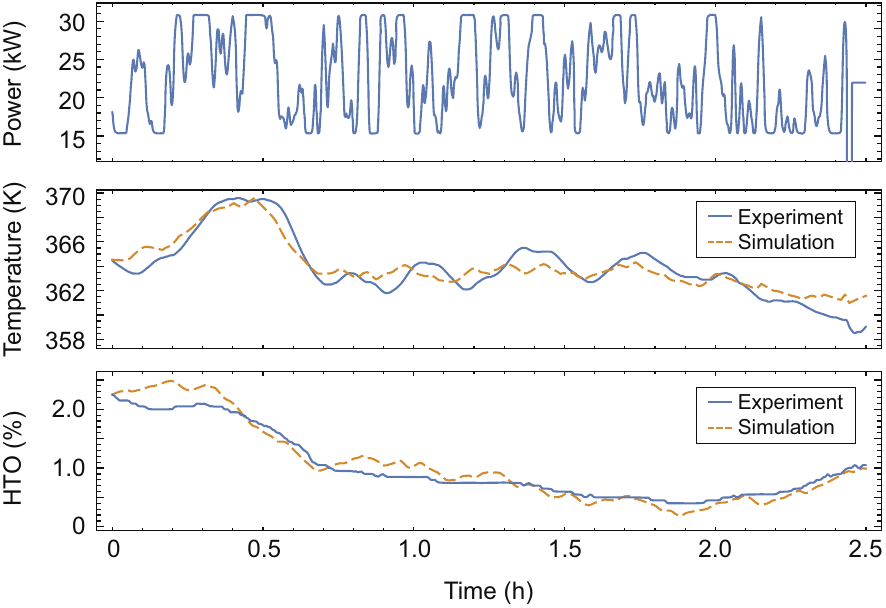}
  \caption{Experimental temperature and HTO impurity crossover data compared with simulations of the dynamic models presented in Section \ref{sec:constraint}.}
  \label{fig:validation}
\end{figure}

The observed temperature, HTO impurity, and simulation results are compared in Fig. \ref{fig:validation}. The parameters of the models are estimated using our previous work \cite{qiu2022online}. We can see that the simulation result fits the experimental data decently. Although our previous works \cite{qi2021pressure} and \cite{qi2022thermal} provide high-order models of the multiphysics effects of the electrolyzer with higher precision, this may further complicate the scheduling problem and make it unsolvable. Hence, we take the models in Sections \ref{sec:temp} and \ref{sec:hto} as acceptable.


\section*{Acknowledgement}

Financial support from National Key Research and Development Program of China (2021YFB4000500), National Natural Science Foundation of China (51907099 and 51907097) is gratefully acknowledged.

\section*{Declaration of Interest}

None.

\section*{Data Availability}

The data related to this work are available upon request.

\section*{References}

\bibliographystyle{elsarticle-num}
\bibliography{AWE-SCHD}

\end{document}